%% file: paper.tex
\documentclass[12pt,reqno,a4paper]{amsart}

\advance\voffset    by -1  cm
\advance\hoffset    by -1.5cm
\advance\textwidth  by  3  cm
\advance\textheight by  1  cm

\input macros

\usepackage{epsfig} 

\begin{document}

\bibliographystyle{plain}

\title{Nonintersecting lattice paths
on the cylinder}

\begin{abstract}
\input abstract
\end{abstract}

\author{Markus Fulmek}
\address{Institut f\"ur Mathematik, Universit\"at Wien\\
Strudlhofgasse 4, A-1090 Wien, Austria}
\email{{\tt Markus.Fulmek@Univie.Ac.At}\newline\leavevmode\indent
{\it WWW}: {\tt http://www.mat.univie.ac.at/\~{}mfulmek}
}

\date{\today}
\thanks{Research partially supported by European Commission's
	IHRP Programme,
	grant HPRN--CT--2001--00272, ``Algebraic Combinatorics in
	Europe''.}
\maketitle

\section{Introduction}
\label{intro}
\input intro

\section{Basic Definitions and Presentation of known Formulas}
\label{basic}
\input basic

\section{Results and Proofs}
\label{results}
\input results

\section{Applications}
\label{applics}
\input applics

\input applics2


\bibliography{paper}

\end{document}

%% file: macros.tex
\def\bit{\begin{itemize}}
\def\eit{\end{itemize}}

\def\wexp{{\mathcal W}}
\def\lgt{Lindstr\"om--Gessel--Viennot}

\def\pas#1{\left(#1\right)}

\def\image{\operatorname{im}\!}

\def\bit{\begin{itemize}}
\def\eit{\end{itemize}}

\def\modulo#1{\;\pas{#1}}
\def\modeq#1#2#3{#1\equiv #2\;\pas{#3}}

\def\N{{\mathbb N}}

\def\R{{\mathbb R}}
\def\Z{{\mathbb Z}}

\def\i{{\mathbf i}}
\def\d{{\mathrm d}}
\def\e{{\mathbf e}}

\newcommand{\SG}{{\mathcal S}}
\newcommand{\defeq}{:=}                 
\newcommand{\fpath}{{\mathbf P}}       
\newcommand{\GF}[1]{\mathbf{GF}\of{#1}}
\newcommand{\sgn}{\operatorname{sgn}}                   
\newcommand{\id}{{\mathrm{id}}}             

\newcommand{\thmref}[1]{Theorem~\ref{#1}}

\newcommand{\figref}[1]{Figure~\ref{#1}}

\newcommand{\setof}[1]{\left\{#1\right\}}
\newcommand{\of}[1]{\left(#1\right)}

\newcommand{\firstcolumn}[1]{
\quad \vdots\quad }
\newcommand{\firstrow}[1]{
\quad \vdots\quad }
\newcommand{\brkof}[1]{\left[#1\right]}

\newcommand{\sgnof}[1]{\operatorname{sgn}\of{#1}}

\newcommand{\edge}[2]{\left[#1\rightarrow #2\right]}

\newcommand{\setpathof}[2]{{\mathcal P}\of{#1,#2}}
\newcommand{\hatsetpathof}[2]{\hat{{\mathcal P}}\of{#1,#2}}
\newcommand{\setnipathof}[2]{{\mathcal P^+}\of{#1,#2}}

\newtheorem{thm}{Theorem}
\newtheorem{prop}{Proposition}
\newtheorem{lem}[thm]{Lemma}
\newtheorem{cor}[thm]{Corollary}

\newtheorem{rem}[thm]{Remark}

%% file: abstract.tex
We show how a formula concerning ``vicious walkers'' (which basically are
nonintersecting lattice paths) on the cylinder
given by P.J.~Forrester can be proved and generalized
by using the \lgt\ method, after having things set up in the
right way. We apply the corresponding results to the
(thermodynamic limit of the) free energy of the ``lock step
model of vicious walkers'', thus completing (and in one instance correcting) the work of Forrester .
Moreover, we also show how a related formula given
by I.~Gessel and C.~Krattenthaler can be obtained from the same
``point of view''.

%% file: intro.tex
In this paper, we consider interesting formulas concerning
nonintersecting lattice paths on the cylinder, and show how the
well--known \lgt\ method provides a common and
quite simple framework for proving them.

We also consider
the asymptotic behaviour of these formulas (i.e., we determine
the thermodynamic limit of the free energy of the ``lock
step model of vicious walkers'') and correct a small
error in the respective formula
\cite[(2.33)]{forrester:vicious-walkers} given by Forrester.

\subsection{Forrester's formula}
In his paper \cite{forrester:vicious-walkers}, Forrester considered
the generating function of certain ``vicious walkers''
\cite[Theorem 2.1]{forrester:vicious-walkers}. The model of
vicious walkers was originally introduced by Fisher \cite{fisher:walks}.
In combinatorial terms, Forrester's formula 
simply gives the enumeration of nonintersecting
lattice paths on a cylindric lattice, expressed as a determinant of certain sums,
but only for the case of an {\em odd\/} number of nonintersecting lattice paths.
Forrester proved this formula using a recurrence relation.

\subsection{Simple framework for Forrester's formula: \lgt}
In our paper, we shall show how the \lgt\ 
framework \cite{lindstroem,gessel-viennot:det} for
directed graphs
can be effortlessly adapted to the cylindric lattice $\Z\times\Z_M$.
From this point of view, Forrester's formula literally is ``easily seen''.

Moreover, it is almost immediate
that in this
setting the appropriate generalization of Forrester's formula also holds for
an {\em even\/} number of ``vicious walkers''. However, in its ``raw''
form, the respective formula
contains summands with negative sign, and hence
is not very useful for enumeration purposes. We overcome this
disadvantage by
appropriately modifying the weights in the respective generating function
(see Theorem~\ref{thm:main2}). 

As applications, we give enumeration formulas for the case of $r$ {\em equidistant\/}
vicious walkers.
While for an odd number $r$, this formula
is already contained in \cite[(2.28)]{forrester:vicious-walkers}, the formula
for even $r$ seems to be new. Thus we are able to complete
Forrester's work; in particular, we can now also determine
the 
asymptotics for even $r$.
Finally, we indicate another proof (basically amounting to coefficient
extraction in Forrester's formula) of a formula given by
Gessel and Krattenthaler \cite{gessel-kratt:cylindric-partitions}.

\subsection{Organization of this paper} 
This paper is organized as follows:

\bit
\item In Section~\ref{basic}, we present the basic definitions and recall
the \lgt\ method.
\item In Section~\ref{results}, we explain how the \lgt\ me\-thod applies to a cylindric lattice.
\item In Section~\ref{applics}, we derive explicit
	enumeration formulas for equidistant vicious walkers.
	(These are related to a formula obtained by Grabiner
	\cite[(33)]{grabiner:random-walks2003}: Our
	formulas \eqref{eq:eqd-odd} and \eqref{eq:eqd-even} could
	be obtained by an appropriate summation of Grabiner«s
	formula). Moreover, we give the corresponding
	asymptotic formulas for the number of paths tending to
	infinity (and correct a small error in the asymptotic formula
	\cite[(2.33)]{forrester:vicious-walkers} given by Forrester).
	Finally, we indicate
	how Forrester's formula and our generalization is related
	to the
	main theorem of Gessel and Krattenthaler
	\cite[Proposition 1, Equation (3.5)]{gessel-kratt:cylindric-partitions}.
\eit

%% file: basic.tex
The main purpose of this paper is to present how the right point of
view almost immediately gives insight in Forrester's formula as well as in
Gessel and Krattenthaler's formula. So we make an effort to give
a careful explanation of this point of view.
\subsection{Nonintersecting paths and generating functions in the
lattice $\Z\times\Z$}

\begin{figure}
\caption{Illustration of lattice paths in $\Z\times\Z$. The picture shows
three lattice paths $p_1$, $p_2$ and $p_3$; from $\pas{3,-2}$ to $\pas{6,11}$,
from $\pas{6,-2}$ to $\pas{11,11}$, and from $\pas{9,-2}$ to $\pas{10,11}$,
respectively. Note that $p_1$ intersects $p_3$ in point $\pas{7,4}$, but
$p_2$ does neither intersect $p_1$ nor $p_3$, since the ``geometric
crossings'' do not correspond to common lattice points.}
\label{fig:zz-lat}
\centerline{\epsfxsize12cm\epsffile{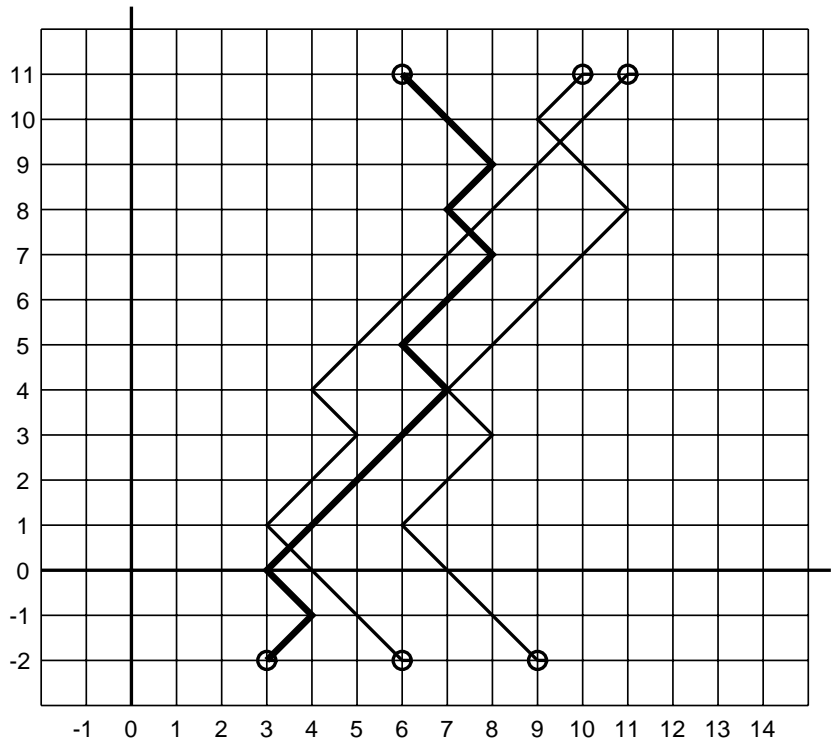}}
\end{figure}

Consider the lattice $\Z\times\Z$, i.e., the directed graphs with
vertex set $\Z\times\Z$ and arcs from $\pas{m,n}$ to $\pas{m-1,n+1}$
(a ``step to the left'') and from $\pas{m,n}$ to $\pas{m+1,n+1}$
(a ``step to the right'') for all $\pas{m,n}\in\Z\times\Z$ (see
\figref{fig:zz-lat}). To all steps to the right, assign weight $1$, and
to all steps to the left, assign weight $x$, i.e., for the edge
$$e=\edge{\pas{m_0,n}}{\pas{m_1,n+1}}$$
we have
\begin{equation}
\label{eq:weightsteps}
w\of{e} =
\begin{cases}
1 & \text{if } m_1 = m_0 + 1, \\
x & \text{if } m_1 = m_0 - 1.
\end{cases}
\end{equation}

A path $p$ of length $N$ is simply a sequence of $N$ adjacent edges
$\pas{e_1,\dots,e_N}$; i.e., for the sequence $\pas{v_i}_{i=0}^N$ of
vertices in the path, we have
$$
e_i = \edge{v_{i-1}}{v_i}
$$
for $i = 1,\dots N$.
The vertices $v_0$ and $v_N$ are called the starting point
and the end point of $p$, respectively.

The weight of a lattice path $p=\pas{e_1,\dots,e_N}$ of length $N$ is
simply defined to be the product of the weights of its edges, i.e.,
\begin{equation}
\label{eq:unsigned-weights}
w\of p = \prod_{i=1}^N w\of{e_i}.
\end{equation}

Two lattice paths $p_1$ and $p_2$ are called {\em intersecting\/}, if they
have a {\em vertex\/} (i.e., a {\em lattice point\/}) in common. A family of lattice paths $p_1,\dots,p_r$
(also called an $r$--tuple of lattice paths) is called {\em nonintersecting\/},
if no two of its paths are intersecting.

See \figref{fig:zz-lat} for an illustration of these simple concepts.

\begin{rem}
\label{rem:evenonly}
Note that ``intersection'' refers only to {\em common lattice
points}: E.g., the ``geometric crossings'' of path $p_1$ and $p_2$ in
\figref{fig:zz-lat} do {\em not\/} constitute intersections
in this sense.

It is clear that two paths starting in lattice points $\pas{m_1,n_1}$ and
$\pas{m_2,n_2}$, respectively, can only be intersecting if
$\pas{m_1-m_2+n_1-n_2}$ is an {\em even\/} number.

While the following considerations and formulas are valid even if this
parity condition is violated, the most interesting case occurs if we
consider the ``even--numbered'' sub--lattice.
\end{rem}

The weight of an arbitrary (not necessarily non--intersecting) $r$--tuple
$\fpath=\pas{p_1,\dots,p_r}$ of lattice paths is simply defined
to be the product of the weights of the single paths, i.e.,
\begin{equation}
w\of{\fpath} = \prod_{i=1}^r w\of{p_i}.
\end{equation}

As usual, by the generating function of some set $A$ of weighted objects
we understand the sum of the weights of the objects, i.e.,
$$
\GF{A} = \sum_{a\in A} w\of{a}.
$$
\subsection{The \lgt\ determinant}
The enumeration of nonintersecting paths in some directed
graph with given starting and
end points is given by the \lgt\
determinant (see \cite[Lemma 1]{lindstroem} or
\cite[Corollary 2]{gessel-viennot:det}). In order to make clear how
this elegant method can be applied to the case of cylindric lattices also,
we state 
this well--known result:
\begin{prop}
\label{prop:lgv-det}
Let $D=\pas{{\mathcal V},{\mathcal A}}$ a directed graph (with vertex
set ${\mathcal V}$ and arc set ${\mathcal A}$), and let $A=\pas{a_1,\dots,a_r}$ and $E=\pas{e_1,\dots,e_r}$ be
two lists of arbitrary vertices in the $D$.
Then we have:
\begin{equation}
\label{eq:lgt-gen}
\det_{1\leq i,j,\leq r}\of{\GF{\setpathof{a_i}{e_j}}} =
\sum_{\pi\in\SG_r}\sgnof\pi
	\GF{\setnipathof{A}{E_{\pi}}}.
\end{equation}
where $\setpathof{a}{e}$ denotes the set of {\em all\/} paths
starting at $a$ and ending at $e$, and $\setnipathof{A}{E_{\pi}}$
denotes the set of all $r$--tuples of {\em nonintersecting\/} paths, where path $i$ starts at $a_i$ and ends at $e_{\pi\of i}$.
\end{prop}

\begin{rem}
\label{rem:lgt}
The ``usual application'' of
Proposition~\ref{prop:lgv-det} 
contains the additional assumption, that
for $1\leq i < j \leq r$ and $1\leq k < l \leq r$, {\em any\/} path from
$a_i$ to $e_l$ must intersect {\em any\/} path from $a_j$ to $e_k$. In this case,
there is only one 
summand on the right--hand side of
\eqref{eq:lgt-gen}, namely $\GF{\setnipathof{A}{E}}$, which
corresponds to the identity permutation.
\end{rem}

\subsection{Paths and generating functions in the lattice $\Z_M\times\Z$}
\begin{figure}
\caption{Illustration  of the $M$--cylinder for $M=12$. The picture shows a
	lattice path $p$ of length $4$ and weight $x$,
	starting in $\pas{0,0}$ and ending in $\pas{2,4}$.}
\label{fig:cyl-lat-path}
\centerline{\epsfxsize6cm\epsffile{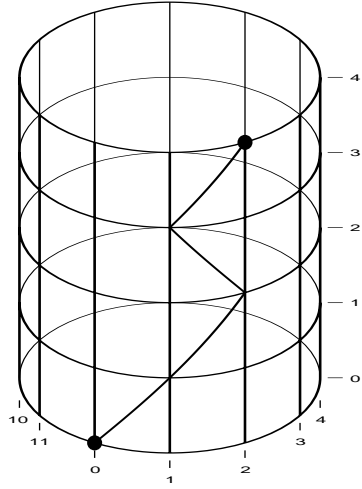}}
\end{figure}

For some (arbitrary, but fixed) integer $M>1$, consider the mapping
\begin{align}
\wexp &: \Z\times\Z\rightarrow\R^3, \notag\\
\wexp\of{m,n}& =\pas{\cos\frac{2\pi m}{M}, \sin\frac{2\pi m}{M}, n}.
	\label{eq:mapcylinder}
\end{align}
Note that the mapping $\wexp$ simply ``wraps'' the lattice $\Z\times\Z$
``around the cylinder''. More precisely,
view the image $\image\,\wexp$ as a ``cylindric lattice''; i.e., a lattice
with vertex set $\setof{\wexp\of{q}:q\in\Z_M\times\Z}$,
and edges leading from points
\begin{align*}
\edge{\wexp\of{m,n}}{\wexp\of{m+1,n+1}} & \text{ (a ``counter--clockwise'' step),}\\
\edge{\wexp\of{m,n}}{\wexp\of{m-1,n+1}} & \text{ (a ``clockwise'' step).}
\end{align*}
We shall call this cylindric lattice the $M$--cylinder.
(\figref{fig:cyl-lat-path} illustrates this simple concept.)

Clearly, a lattice path $p$ in the $M$--cylinder can be viewed as the 
image of an ``ordinary'' lattice path in $\Z\times\Z$ under the mapping $\wexp$.
Each path $p$ in the $M$--cylinder inherits the weight from a corresponding
path in $\Z\times\Z$, i.e., if $p=\wexp\of{\hat{p}}$, we set
$w\of p \defeq w\of{\hat{p}}$. (Note that this
is well--defined.)

\begin{figure}
\caption{Illustration of intersecting lattice paths on the $M$--cylinder for $M=12$. The right picture shows representatives
of the preimages of these paths (drawn with thick lines)
in the lattice
$\Z\times\Z$ under the mapping defined in \eqref{eq:mapcylinder}.
The whole
preimage consists of an infinite family of horizontally
translated paths, indicated by 
thin lines in the picture.}
\label{fig:int-paths}
\centerline{\epsfxsize13cm\epsffile{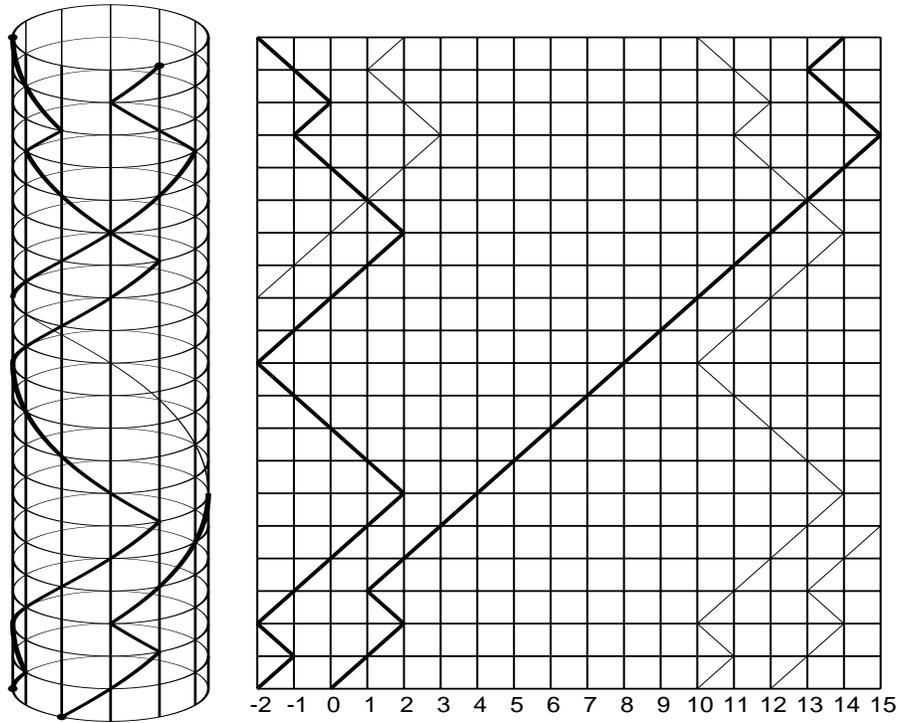}}
\end{figure}

A lattice path may ``wind
around the cylinder several times'', in either positive or negative direction,
before reaching its end point:
The preimage
$\wexp^{-1}\of{\setpathof{a}{e}}$ of the set of
lattice paths in the $M$--cylinder,
which start at $\pas{a,0}$
and end in $\pas{e,N}$, consists of
lattice paths in $\Z\times\Z$, which start at $\pas{a+k\cdot M,0}$
and end in $\pas{e+\pas{k+o}\cdot M,N}$ for $k,o\in\Z$, i.e.,
\begin{equation*}
\wexp^{-1}\of{\setpathof{a}{e}} =
\bigcup_{o\in\Z}\pas{
	\bigcup_{k\in\Z}\hatsetpathof{a+k\cdot M}{e+\pas{k+o}\cdot M}
}.
\end{equation*}
We shall call the number $o$ the
{\em offset of the
endpoint\/} of the path $p$. (See \figref{fig:int-paths} for this
concept.)



In the following, we shall restrict ourselves to lattice paths starting
at $\pas{a,0}$ and ending at $\pas{e,N}$, where $N$ is some (arbitrary
but fixed) integer. Note that in $\Z\times\Z$, a lattice path
starting at $\pas{a,0}$ and ending at $\pas{e,N}$ exists
if and only if
\begin{equation}
N-2k=e-a
\label{eq:condplain}
\end{equation}
for some $k\in\Z$ with $0\leq k\leq N$ (here, $k$ denotes the number of steps to the left).
For lattice paths in the $M$--cylinder, this condition is changed to
\begin{equation}
N-2k=\pas{e+o\cdot M}-a
\label{eq:condcylinder}
\end{equation}
for arbitrary $o\in\Z$ and some $k\in\Z$  with $0\leq k\leq N$
(here, $o$ denotes the {\em offset\/} of the
endpoint, and $k$ denotes the number of clockwise steps).
 
So it is easy to see that
the generating function $q\of{M,N,a,e}$ of all lattice paths
in the $M$--cylinder, which start at $\pas{a,0}$ and end in
$\pas{e,N}$, is given by
\begin{equation}
\label{eq:gf_qx}
q\of{M,N,a,e;\;x}
=
\sum_{\shortstack{
	$\scriptstyle o\in\Z$\\
	$\scriptstyle \modeq{N-e-o\cdot M+a}{0}{2}$
	}
} \binom{N}{\frac{N-e-o\cdot M+a}{2}}x^{\frac{N-e-o\cdot M+a}{2}}.
\end{equation}
Forrester gave an equivalent expression (see
\cite[equation 2.11 and 2.12]{forrester:vicious-walkers}) for \eqref{eq:gf_qx}, which is more elegant insofar as it ``conceals'' the
clumsy definition of the range of  summation in \eqref{eq:gf_qx}:
\begin{align}
\label{eq:gf_qx-forrester}
q\of{M,N,a,e;\;x}
&=
\frac{1}{M}
\sum_{l=0}^{M-1}
	\e^{\frac{-2\pi\i\pas{e-a}l}{M}}
		\pas{x\e^{\frac{-2\pi\i l}{M}}+\e^{\frac{2\pi\i l}{M}}}^N
\\ 
& =
\sum_{k=0}^N \binom{N}{k}x^k
	\frac{1}{M}
	\sum_{l=0}^{M-1} \pas{\e^{\frac{2\pi\i}{M}\pas{N-2k-e+a}}}^l, \notag
\end{align}
where $\i$ denotes the imaginary unit.
\eqref{eq:gf_qx-forrester} is equal to \eqref{eq:gf_qx}, since
we have
\begin{equation}
\label{eq:unit-roots-powersum}
\sum_{l=0}^{M-1} \pas{\e^{\frac{2\pi\i}{M} m}}^l =
\begin{cases}
	M &\text{ if } m\equiv 0\modulo{M}, \\
	0 &\text{ else.}
\end{cases}
\end{equation}

%% file: results.tex
It is an obvious observation that the \lgt\ method for nonintersecting paths in a directed graph, as described in
Proposition~\ref{prop:lgv-det}, applies to the $M$--cylinder.
However, unlike the ``usual case'' outlined
in Remark~\ref{rem:lgt}, 
there appear terms corresponding to {\em other\/} permutations
than the identity permutation; which turn out to be
{\em cyclic\/} permutations.
These permutations may
have {\em negative sign\/} if the number of paths, $r$, is even.
So if we are given starting points
$\pas{\pas{a_1,0},\dots\pas{a_r,0}}$ and end points $\pas{\pas{e_1,N},\dots\pas{e_r,N}}$,
we define the {\em signed weight\/} of a family $\fpath=\pas{P_1,\dots,P_r}$
of lattice paths, where $P_i$ starts in $\pas{a_i,0}$ and ends
in $\pas{e_{\pi\of i},N}$ for some permutation $\pi\in\SG_r$, as
\begin{equation}
w\of{\pi,\fpath} = \sgn\pi\prod_{j=1}^r w\of{P_j}.
\label{eq:signed-weights2}
\end{equation}

\subsection{A simple generalization of Forrester's formula}
Given this ``signed weight'' for families of nonintersecting lattice paths,
we may derive immediately the following generalization of Forrester's formula.
\begin{thm}
\label{thm:main}
The generating function {\em with signed weights\/} (according to
\eqref{eq:signed-weights2})
of all $r$--tuples of non--intersecting lattice
paths in the $M$--cylinder, starting at the points
$$\pas{\pas{a_1,0},\dots\pas{a_r,0}}$$
and ending in any permutation of the points
$$\pas{\pas{e_1,N},\dots\pas{e_r,N}},$$
with $0\leq a_1< a_2<\dots<a_r<M$ and $0\leq e_1<e_2 <\dots < e_r<M$,
where for all $1\leq i,j\leq r$ we have $\modeq{N-e_i+a_j}{0}{2}$,
is given by
\begin{equation}
\label{eq:det-gfx}
\det\of{q\of{M,N,a_i,e_j;x}}_{i,j=1}^{r}.
\end{equation}
Moreover, we have the following expansion for the above
determinant:
\begin{equation}
\label{eq:det-gfx-expansion}
\det\of{q\of{M,N,a_i,e_j;x}}_{i,j=1}^{r}=
\sum_{i=0}^{r-1}\sgnof{\mu^i}
	\GF{\setnipathof{A}{E_{\mu^i}}},
\end{equation}
where $\mu$ denotes the permutation mapping $1$ to $2$,
$2$ to $3$, and so on; i.e., $\mu=\pas{1,2,\dots,r}$ in
cycle notation.
\end{thm}

\begin{proof}
The assertion of \eqref{eq:det-gfx} is an immediate consequence
of Proposition~\ref{prop:lgv-det}.


For the assertion of \eqref{eq:det-gfx-expansion},
consider the {\em preimage\/} of any nonintersecting
$r$--tuple of lattice paths with respect to $\wexp$: It appears as a periodic
configuration of {\em infinitely many\/} nonintersecting lattice paths in
$\Z\times\Z$, such that each point $\pas{a_i+s\cdot M,0}$ and each point
$\pas{e_i+s\cdot M,N}$ (for $i=1,\dots,r$ and $s\in\Z$) appears as starting
point and as end point, respectively, of some path (see \figref{fig:preimage}
for an illustration).

\begin{figure}
\caption{Illustration of the {\em preimage\/}
with respect to $\wexp$ (right picture)
of a nonintersecting
pair of lattice paths in the $6$--cylinder (left picture).}
\label{fig:preimage}
\centerline{\epsfxsize16cm\epsffile{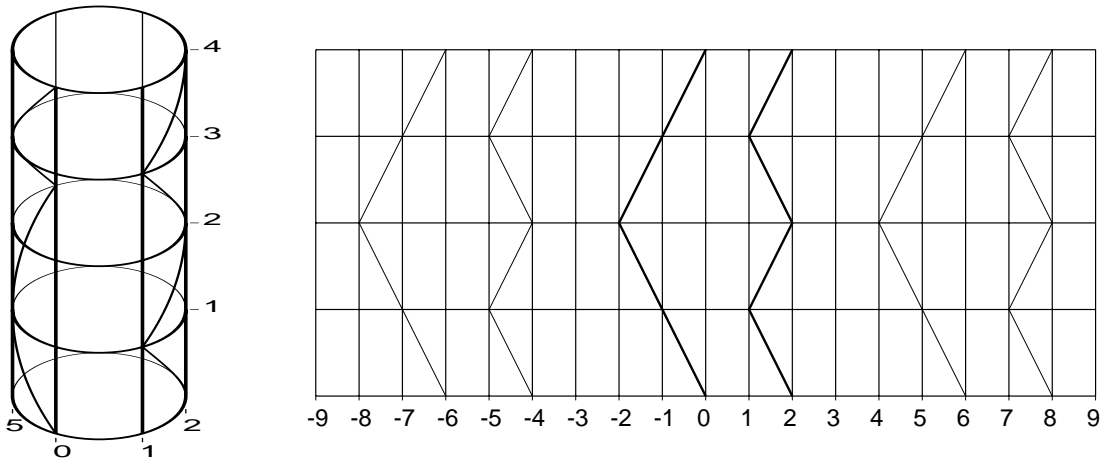}}
\end{figure}

It is obvious that
such a configuration can only correspond to some ``shift of the endpoints'' in the following sense: Consider the ``canonical'' starting points $\pas{a_i,0}$, and
label the ``canonical'' endpoints $\pas{e_i,N}$ with the numbers $1,2,\dots,r$.
Label the other possible endpoints
from left to right with the integers in a consistent way (i.e., $\pas{e_r-M,N}$ gets label
0, $\pas{e_{r-1}-M,N}$ gets label $-1$; $\pas{e_1+M,N}$ gets label
$r+1$, and so on). Then for any nonintersecting $r$--tuple of lattice paths there
is some fixed integer $p$, such that for $1\leq i\leq r$, the path starting
at $\pas{a_i,0}$ ends in the $(i+p)$--th endpoint in this labeling. This ``shift of the endpoints'' clearly
corresponds to a cyclic permutation $\mu^p$, where $\mu=\pas{1,2,\dots,r}$.
\end{proof}

%
\begin{rem}
Note that under the assumptions of Theorem~\ref{thm:main}, the
case $M$ odd admits only one possible permutation of endpoints,
namely the identity permutation, since all such permutations
must be of the form $\pas{\mu^r}^{2k}=\id$, 
according to condition
\eqref{eq:condcylinder}.

While the following considerations and formulas are valid also for odd $M$,
the most interesting cases occur if we assume
\bit
\item $M\equiv N \equiv 0\pas{2}$,
\item $a_i-a_j\equiv e_i-e_j\equiv 0\pas{2}
	\;\forall 1\leq i,j\leq r$.
\eit
(See also Remark~\ref{rem:evenonly}.)
\end{rem}

The determinantal expression \eqref{eq:det-gfx}
was derived
recursively by Forrester only for odd $r$ (see \cite[Theorem 2.1, equation 2.10]{forrester:vicious-walkers}). In \thmref{thm:main}, we easily extended
it to the case of even $r$ (thus answering the respective question
posed in \cite{forrester:vicious-walkers}), just by adopting 
the right ``point of view''
(i.e., the \lgt\ framework with ``signed weights'').

\subsection{A more sophisticated generalization of Forrester's formula}
Note that for odd $r$, all the (cyclic) permutations in
\eqref{eq:det-gfx-expansion}
have positive sign. 
If $r$ is even, however, also negative terms appear
in the generating
function \eqref{eq:det-gfx}: This is a bit of a
nuisance, for we cannot simply set $x\equiv1$ in order
to obtain an {\em enumeration formula\/}.

An easy way out of this
difficulty is to modify the definition of the weight of a single path $p$,
so that
its ``offset of endpoint'', $o$, is taken into account via a multiplicative
factor of $y^o$, i.e.,
we replace definition
\eqref{eq:unsigned-weights} by
\begin{equation}
\label{eq:y-weights}
w_y\of p = y^o\prod_{i=1}^N w\of{e_i}.
\end{equation}
This amounts to the following
modification of \eqref{eq:gf_qx}:
\begin{equation}
\label{eq:gf_qxy}
q\of{M,N,a,e;\;x,y}
=
\sum_{\shortstack{
	$\scriptstyle o\in\Z$\\
	$\scriptstyle \modeq{N-e-o\cdot M+a}{0}{2}$
	}
} \binom{N}{\frac{N-e-o\cdot M+a}{2}}y^o\cdot x^{\frac{N-e-o\cdot M+a}{2}}.
\end{equation}
Now, the proof of Theorem~\ref{thm:main} (with slight and obvious
modifications) immediately yields the following Lemma:
\begin{lem}
\label{lem:main}
The generating function {\em with weights according to \eqref{eq:y-weights}\/}
of all $r$--tuples of non--intersecting lattice
paths in the $M$--cylinder, starting at the points
$$\pas{\pas{a_1,0},\dots\pas{a_r,0}}$$
and ending in some permutation of the points
$$\pas{\pas{e_1,N},\dots\pas{e_r,N}},$$
with $0\leq a_1< a_2<\dots<a_r<M$ and $0\leq e_1<e_2 <\dots < e_r<M$,
where for all $1\leq i,j\leq r$ we have $\modeq{N-e_i+a_j}{0}{2}$,
is given by
\begin{equation}
\label{eq:det-gfxy}
\det\of{q\of{M,N,a_i,e_j;\;x,y}}_{i,j=1}^{r}.
\end{equation}
\end{lem}
A simple argument now leads to the desired formula without unwanted negative
signs:
\begin{thm}
\label{thm:main2}
The generating function {\em with unsigned weights\/} (i.e., according to \eqref{eq:unsigned-weights})
of all $r$--tuples of non--intersecting lattice
paths in the $M$--cylinder, starting at the points
$$\pas{\pas{a_1,0},\dots\pas{a_r,0}}$$
and ending in some permutation of the points
$$\pas{\pas{e_1,N},\dots\pas{e_r,N}},$$
with $0\leq a_1< a_2<\dots<a_r<M$ and $0\leq e_1<e_2 <\dots < e_r<M$,
where for all $1\leq i,j\leq r$ we have $\modeq{N-e_i+a_j}{0}{2}$,
is given by
\begin{equation}
\label{eq:det-gfx-unsigned}
\det\of{q\of{M,N,a_i,e_j;\;x,\pas{-1}^{r-1}}}_{i,j=1}^{r}.
\end{equation}
\end{thm}
\begin{proof}
For odd $r$, we clearly have
$$
q\of{M,N,a,e;\;x,1} = q\of{M,N,a,e;\;x},
$$
whence \eqref{eq:det-gfx-unsigned} simply amounts to the assertion of
Theorem~\ref{thm:main}.

For even $r$, observe that (due to \eqref{eq:det-gfx-expansion})
the sign of the summands in
\eqref{eq:det-gfxy}
equals $\pas{-1}^{n}$, where $n$ is the number of paths with {\em odd\/}
offset of endpoint in the corresponding $r$--tuple of paths. So setting
$y=-1$ in \eqref{eq:det-gfxy} properly cancels all the negative signs.
\end{proof}

\subsection{Enumeration formulas involving trigonometric functions}
It is possible to rewrite the generating function
$q\of{M,N,a,e;x,y}$
in a way similar to \eqref{eq:gf_qx-forrester}.

\begin{cor}
\label{cor:main}
For $M>0$, we have:
\begin{equation}
\label{eq:cor4}
q\of{M,N,a,e;x,y}
= 
\frac{y^{\frac{a-e}{M}}}{M}
\sum_{l=0}^{M-1}
	\e^{\frac{-2\pi\i\pas{e-a}l}{M}}
		\pas{
			x y^{-\frac{1}{M}}\e^{\frac{-2\pi\i l}{M}}+
			y^{\frac{1}{M}}\e^{\frac{2\pi\i l}{M}}
		}^N.
\end{equation}
Moreover, the generating function \eqref{eq:det-gfx-unsigned} of
Theorem~\ref{thm:main2} is equivalently given by
\begin{equation}
\label{eq:det-gfxr}
M^{-r}
\det\Biggl(
	\e^{\frac{\pas{r-1}\pi\i\pas{a_i-e_j}}{M}}
	\sum_{l=0}^{M-1}
		\e^{\frac{-2\pi\i\pas{e_j-a_i}l}{M}}
			\pas{
				x\e^{\frac{-2\pi\i\pas{l+\frac{r-1}{2}}}{M}}+
				\e^{\frac{2\pi\i\pas{l+\frac{r-1}{2}}}{M}}
			}^N
\Biggr)_{i,j=1}^{r}.
\end{equation}
\end{cor}
\begin{proof}
Equation \eqref{eq:cor4} follows from the same type of computation as in
\eqref{eq:gf_qx-forrester}.
\input derivation0
Now set $y=\pas{-1}^{r-1}=\e^{\pas{r-1}\pi\i}$ in \eqref{eq:cor4}, and insert the result into
\eqref{eq:det-gfx-unsigned}: This immediately yields
\eqref{eq:det-gfxr}.
\end{proof}

If the starting points $a_i$ 
are {\em equidistant\/},
we can simplify the corresponding expressions
even
further by a little trick, extending the computation carried
out by Forrester
\cite[equation (3.2)]{forrester:vicious-walkers} to the
case of even $r$:
\begin{cor}
\label{cor:eq-start}
Consider the case of equidistant starting points, i.e., let
$a_i=\pas{i-1}\cdot\nu$ and 
$M=r\cdot\nu$
in\eqref{eq:det-gfxy} for some fixed $\nu\in\N$. In this case,
the generating function \eqref{eq:det-gfx-unsigned} of
Theorem~\ref{thm:main2} is given by
\begin{multline}
\label{eq:eqdiststart-even}
\i^{-\frac{(r)(r-1)}{2}}\pas{\nu\sqrt{r}}^{-r}\\\times
\det\of{\sum_{a=0}^{\nu-1}
		\e^{\frac{\pi\i\pas{N-e_{j+1}\pas{2\pas{i+a r}+1}}}{r\nu}}
			\pas{
				x\e^{\frac{-2\pi\i\pas{i+a r+1}}{r\nu}}+
				\e^{\frac{2\pi\i\pas{i+a r}}{r\nu}}
			}^N}_{i,j=0}^{r-1}
\end{multline}
if $r$ is even, and by
\begin{equation}
\label{eq:eqdiststart-odd}
\i^{-\frac{(r+2)(r-1)}{2}}\pas{\nu\sqrt{r}}^{-r}
\det\of{\sum_{a=0}^{\nu-1}
		\e^{\frac{\pi\i\pas{-e_{j+1}\pas{2\pas{i+a r}}}}{r\nu}}
			\pas{
				x\e^{\frac{-2\pi\i\pas{i+a r}}{r\nu}}+
				\e^{\frac{2\pi\i\pas{i+a r}}{r\nu}}
			}^N}_{i,j=0}^{r-1}
\end{equation}
if $r$ is odd. (Note that --- as a matter of convenience ---
row and column indices range
from $0$ to $\pas{r-1}$ here.)
\end{cor}
\begin{proof}
Note that
$\brkof{
		\pas{\e^{-\frac{2\pi\i}{r}i}}^j
		r^{-1/2}
}_{i,j=0}^{r-1}$
is a unitary matrix. This fact, together with the well--known
formula for Vandermonde determinants, yields the determinant
evaluations
\begin{equation}
\label{eq:some-determinants}
\det\of{\e^{-\frac{2\pi\i}{r}i\cdot j}}_{i,j=0}^{r-1}=
	r^{\frac{r}{2}}\i^{\frac{(r+2)(r-1)}{2}},
\text{ and }
\det\of{\e^{-\frac{2\pi\i}{r}\pas{i+1/2}\cdot j}}_{i,j=0}^{r-1}=
	r^{\frac{r}{2}}\i^{\frac{(r)(r-1)}{2}}.
\end{equation}
The assertions follow by a simple computation, which we shall
show for even $r$ only (the case $r$ odd being completely analogous; see Forrester
\cite[equation (3.2)]{forrester:vicious-walkers}).
Set $y=\pas{-1}=\e^{\pi\i}$ in \eqref{eq:cor4} and insert in
\eqref{eq:det-gfx-unsigned}.
Now multiply this with the second determinant
from \eqref{eq:some-determinants};
i.e., consider the determinant of the product of the
$r\times r$--matrices
$$
\brkof{\e^{-\frac{2\pi\i}{r}\pas{k+1/2}\cdot m}}
\times
\brkof{
\frac{\e^{\pi\i\frac{N-e_{j+1}+i\nu}{r\nu}}}{r\nu}
\sum_{l=0}^{r\nu-1}
	\e^{\frac{-2\pi\i\pas{e_{j+1}-i\nu}l}{r\nu}}
		\pas{x \e^{\frac{-2\pi\i \pas{l+1}}{r\nu}}+\e^{\frac{2\pi\i l}{r\nu}}}^N
},
$$
where the row and column indices $i$ and $j$ range from
$0$ to $(r-1)$.
The $(i,j)$--entry of this matrix product is given by
\begin{align*}
a_{i,j}&=
\sum_{n=0}^{r-1}\pas{
	\e^{-\frac{2\pi\i}{r}\pas{i+1/2} n}
	\frac{\e^{\pi\i\frac{N-e_{j+1}+n\nu}{r\nu}}}{r\nu}
	\sum_{l=0}^{r\nu-1}
	\e^{\frac{-2\pi\i\pas{e_{j+1}-n\nu}l}{r\nu}}
		\pas{x \e^{\frac{-2\pi\i \pas{l+1}}{r\nu}}+
			\e^{\frac{2\pi\i l}{r\nu}}}^N
}\\
&=
\frac{1}{r\nu}
\sum_{l=0}^{r\nu-1}
	\pas{x \e^{\frac{-2\pi\i \pas{l+1}}{r\nu}}+
		\e^{\frac{2\pi\i l}{r\nu}}}^N
	\e^{\frac{\pi\i}{r\nu}
		\pas{N-e_{j+1}\pas{2l+1}}
	}
	\sum_{n=0}^{r-1}
		\e^{-\frac{2\pi\i}{r}\pas{i-l} n}\\
&=
\frac{1}{\nu}
\sum_{a=0}^{\nu-1}
	\pas{x \e^{\frac{-2\pi\i \pas{i+a r+1}}{r\nu}}+
		\e^{\frac{2\pi\i \pas{i+a r}}{r\nu}}}^N
	\e^{\frac{\pi\i}{r\nu}
		\pas{N-e_{j+1}\pas{2\pas{i+a r}+1}},
	}
\end{align*}
which immediately gives \eqref{eq:eqdiststart-even}. (The proof of
\eqref{eq:eqdiststart-odd} involves multiplication with
the first determinant from \eqref{eq:some-determinants}.)
\end{proof}

%% file: derivation0.tex
\begin{align*}
&\phantom{=}
\frac{y^{\frac{a-e}{M}}}{M}
\sum_{l=0}^{M-1}
	\e^{\frac{-2\pi\i\pas{e-a}l}{M}}
		\pas{
			x y^{\frac{-1}{M}}\e^{\frac{-2\pi\i l}{M}}+
			y^{\frac{1}{M}}\e^{\frac{2\pi\i l}{M}}}^N
  \\
  &=
\frac{y^{\frac{a-e}{M}}}{M}
\sum_{l=0}^{M-1}
	\e^{\frac{-2\pi\i\pas{e-a}l}{M}}
	\sum_{k=0}^N\binom{N}{k}
		x^k y^{\frac{N-2k}{M}}\e^{\frac{2\pi\i}{M}\pas{N-2k}l}
  \\
  &=
\sum_{k=0}^{N}\binom{N}{k}x^k y^{\frac{N-e+a-2k}{M}}
	\frac{1}{M}
	\sum_{l=0}^{M-1}\e^{\frac{2\pi\i}{M}\pas{N-e+a-2k}l}.
\end{align*}
Use \eqref{eq:unit-roots-powersum} to see that
\eqref{eq:cor4}
equals \eqref{eq:gf_qxy} (set 
$k=\frac{N-e+a-o\cdot M}{2}$, or, equivalently, 
$o=\frac{N-e+a-2k}{M}$).

%% file: applics.tex
\subsection{An enumeration formula for a special case}
Of particular interest is the enumeration formula for the case
$M = r\nu$ with equidistant starting points and end points, $a_i=e_i=\pas{i-1}\nu$:
This simply amounts to setting $M = r\nu$, $a_i=e_i=\pas{i-1}\nu$
and $x=1$ in \eqref{eq:det-gfx-unsigned}.
Since we have the obvious relations
\begin{equation}
\label{eq:qsym}
q\of{M,N,a,e;x,y}=x^N q\of{M,N,e,a;1/x,1/y}
\end{equation}
and
\begin{equation}
\label{eq:qtrans}
q\of{r \nu,N,i \nu,j \nu;x,y} =
	y\cdot q\of{r\ \nu,N,i \nu,(j+r) \nu;x,y},
\end{equation}
we may concentrate on the numbers
$$
	a_d\defeq q\of{r \nu,N,0,d \nu;1,\pas{-1}^{r-1}}
	\text{ for } d=0,\dots,r-1.
$$
So for odd $r$, 
we obtain a circulant matrix, the determinant of
which we can easily
evaluate by the well--known formula
(cf.~\cite[\S51, p.~131]{aitken:det}):
\begin{equation}
\label{eq:circulant}
\det
\begin{pmatrix}
a_0     & a_1    & \dots & a_{r-2} & a_{r-1} \\
a_{r-1} & a_0    & \dots & a_{r-3} & a_{r-2} \\
\vdots  & \vdots &       & \vdots  & \vdots \\
a_1     & a_2    & \dots & a_{r-1} & a_0
\end{pmatrix}
=\prod_{m=0}^{r-1}
	\pas{\sum_{k=0}^{r-1}\pas{\e^{\frac{2m\pi\i}{r}}}^k a_k}.
\end{equation}
%
For even $r$, however, we obtain
a ``skew--symmetric'' circulant matrix (due to 
\eqref{eq:qtrans}), the determinant of
which we can evaluate in much
the same way as \eqref{eq:circulant}. Since this evaluation
appears to be not so well--known, we state and prove it in the
following lemma:
\begin{lem}
\label{lem:skew-circulant}
For arbitrary variables $a_0\dots a_{r-1}$, we have
\begin{equation}
\label{eq:skew-circulant}
\det
\begin{pmatrix}
a_0     & a_1    & \dots & a_{r-2} & a_{r-1} \\
-a_{r-1} & a_0    & \dots & a_{r-3} & a_{r-2} \\
\vdots  & \vdots &       & \vdots  & \vdots \\
-a_1     & -a_2    & \dots & -a_{r-1} & a_0
\end{pmatrix}
=\prod_{m=0}^{r-1}
	\pas{\sum_{k=0}^{r-1}\pas{\e^{\frac{(2 m+1)\pi\i}{r}}}^k a_k}.
\end{equation}
\end{lem}
\begin{proof}
Set $\omega_m\defeq\e^{\frac{\pas{2m+1}\pi\i}{r}}$ and consider
the $r$ vectors
$\vec\omega_m\defeq\pas{\omega_m^0,\omega_m^1,
	\dots,\omega_m^{r-1}}$ for $m=0,\dots,r-1$.
Note that
$\omega_m^r=-1$ and compute the $i$--th component
in the product $A\cdot\vec\omega_m$ (where $A$, of course,
denotes the
matrix in \eqref{eq:skew-circulant}):
\begin{align*}
\pas{A\cdot\vec\omega_m}_i &=
	-a_{r-i}\omega_m^0 -a_{r-i+1}\omega_m^1-\dots
		-a_{r-1}\omega_m^{i-1}+
	a_0\omega_m^{i}+a_1\omega_m^{i+1}+\dots
		a_{r-i-1}\omega_m^{r-1}\\
&= a_0\omega_m^{i}+a_1\omega_m^{i+1}+\dots
		a_{r-i-1}\omega_m^{r-1}+
	a_{r-i}\omega_m^r +a_{r-i+1}\omega_m^{r+1}+\dots
		+a_{r-1}\omega_m^{r+i-1}\\
&= \omega_m^{i}\cdot \pas{\sum_{k=0}^{r-1}a_k\omega_m^k}.
\end{align*}
This shows that $\vec\omega_m$ is an eigenvector of $A$
to the eigenvalue $\pas{\sum_{k=0}^{r-1}a_k\omega_m^k}$, which
proves the assertion.
\end{proof}

\begin{cor}
Denote the number
of all $r$--tuples of non--intersecting lattice
paths in the $\pas{r \nu}$--cylinder, starting at the points
$$\pas{\pas{0,0},\dots\pas{\pas{r-1}\nu,0}}$$
and ending in any permutation of the points
$$\pas{\pas{0,N},\dots\pas{\pas{r-1}\nu,N}},$$
by  $Z\of{N,r,\nu}$. Then we have the following formulas:
\begin{align}
\label{eq:eqd-odd}
Z\of{N,2r-1,\nu}&=
\pas{\frac{2^{N}}{\nu}}^{\!\!2r-1}
\prod_{m=0}^{2r-2}
   \sum_{l=0}^{\nu-1}
       \cos^{N}\of{
       	2\,\pi
			\pas{
				\frac{m}{\nu\pas{2r-1}}
				+
				\frac{l}{\nu}
			}
        },\\
\label{eq:eqd-even}
Z\of{N,2r,\nu}&=
\pas{\frac{2^{N}}{\nu}}^{\!\!2r}
\prod_{m=0}^{2r-1}
   \sum_{l=0}^{\nu-1}
       \cos^{N}\of{
       	2\,\pi
			\pas{
				\frac{m+1/2}{\nu\pas{2r}}
				+
				\frac{l}{\nu}
			}
        }.
\end{align}
\end{cor}

\begin{proof}
Clearly, \eqref{eq:eqd-odd} will follow by simplifying
\eqref{eq:circulant}, and \eqref{eq:eqd-even} will follow by simplifying \eqref{eq:skew-circulant}. We shall give the
corresponding computation for \eqref{eq:eqd-even} only,
the other case is completely analogous.

Straightforward insertion of
$$
	a_d\defeq q\of{r \nu,N,0,d \nu;1,\e^{\pi \i}}
$$
into
\eqref{eq:skew-circulant} gives the following expression:
\begin{equation}
\prod_{m=0}^{2r-1}
	\sum_{k=0}^{2r-1}
   		\pas{
			\pas{\e^{\frac{\pas{2m+1}\pi\i}{2r}}}^k
			\frac{\e^{\frac{\pas{N-k\nu}\pi\i}{2\nu r}}}{2\nu r}
   			\sum_{l=0}^{2\nu r-1}
				\e^{-\frac{k\nu l\pi\i}{\nu r}}
				\pas{
					\e^{-\frac{\pas{l+1}\pi\i}{\nu r}}
					+
					\e^{\frac{l\pi\i}{\nu r}}
				}^N
   		}.
\end{equation}
Now write
$$
\pas{
\e^{-\frac{\pas{l+1}\pi\i}{\nu r}}
+
\e^{\frac{l\pi\i}{\nu r}}
} =
\e^{-\frac{\pi\i}{2\nu r}}2\cos\of{\frac{\pas{l+1/2}\pi}{\nu r}},
$$
pull out appropriate factors, simplify, and interchange
summation; in order to obtain
\begin{equation*}
\pas{\frac{2^{N-1}}{\nu r}}^{2r}
\prod_{m=0}^{2r-1}
	\sum_{l=0}^{2\nu r-1}
		\cos^N\of{\frac{\pas{l+1/2}\pi}{\nu r}}
		\sum_{k=0}^{2r-1}\pas{\e^{\frac{2\pas{m-l}\pi\i}{2r}}}^k.
\end{equation*}
Observe that \eqref{eq:unit-roots-powersum} applies
to the innermost sum,
whence \eqref{eq:eqd-even} follows.  
\end{proof}

\begin{rem}
Equation \eqref{eq:eqd-odd} is basically the same as
Forrester's formula \cite[(2.28)]{forrester:vicious-walkers}.
\end{rem}

\subsection{Free energy}
In his paper, Forrester
considers the dimensionless {\em free energy per unit length\/}
on
a strip--shaped lattice of infinite width and height $N$
(see \cite[(2.30)]{forrester:vicious-walkers}):
\begin{equation}
\label{eq:free-energy}
f_N\of\nu=-\frac{1}{\nu r}\lim_{r\to\infty}\log\of{Z\of{N,r,\nu}}.
\end{equation}
We can apply his considerations now also to the case of
even $r$. According to \eqref{eq:eqd-odd} and \eqref{eq:eqd-even},
respectively, we have
\begin{multline}
-\frac{1}{r\nu}\log\of{Z\of{N,r,\nu}} = \\
-\frac{1}{r\nu}
\Biggl(
	r\pas{N \log 2 - \log \nu}
	+\sum_{m=0}^{r-1}\log\of{
		\sum_{l=0}^{\nu-1}
			\cos^N\of{2\pi\frac{l+\frac{m+\epsilon\of r}{r}}{\nu}}
	}
\Biggr),
\end{multline}
where $\epsilon\of r=1/2$ if $r$ is even, and $\epsilon\of r=0$
if $r$ is odd. In both cases, observe that we have Riemann sums,
which tend to the same integral in the limit:
\begin{equation}
\label{eq:asym-forrester}
-\frac{1}{\nu}
\Biggl(
	\pas{N \log 2 - \log \nu}
	+\int_{0}^{1}\log\of{
		\sum_{k=0}^{\nu-1}
			\cos^N\of{2\pi\frac{k+t}{\nu}}d\,t
	}
\Biggr).
\end{equation}
(This corresponds to Forresters formula
\cite[(2.31)]{forrester:vicious-walkers}.)

Now, following Forrester
\cite[Section 2.4]{forrester:vicious-walkers}, we
consider the {\em free energy per lattice site\/} in the
two--dimensional thermodynamic limit, i.e., the quantitity
$F_\nu\defeq\lim_{N\to\infty}\frac{f_N\of\nu}{N}$.
Of course, the basic idea for evaluating this limit is
``Pull out the dominating term from the sum in
\eqref{eq:asym-forrester}''. However, we must be careful in 
determining this dominating term (there seems to be a
small flaw in
Forresters formula \cite[2.33]{forrester:vicious-walkers}
with respect to this):

For even $\nu$, the dominating term is
\bit
\item $\cos\of{\frac{2\pi t}{\nu}}$
	for $0\leq t\leq \frac12$,
\item $\cos\of{\frac{2\pi\pas{t+\nu-1}}{\nu}}=\cos\of{\frac{2\pi\pas{t-1}}{\nu}}$
	for $\frac12\leq t\leq 1$,
\eit
whence we obtain
\begin{align}
F_\nu & = -\frac{\log 2}{\nu}\notag\\
&-
	\lim_{N\to\infty}\frac{1}{\nu N}
	\int_{0}^{\frac12}\log\of{
		\cos\of{\frac{2\pi t}{\nu}}^N\times
		\sum_{k=0}^{\nu-1}
			\pas{
				\frac{
					\cos\of{2\pi\frac{k+t}{\nu}}
				}{
					\cos\of{\frac{2\pi t}{\nu}}
				}
			}^N
	}\d\,t\notag\\
&-
	\lim_{N\to\infty}\frac{1}{\nu N}
	\int_{\frac12}^{1}\log\of{
		\cos\of{\frac{2\pi\pas{t-1}}{\nu}}^N\times
		\sum_{k=0}^{\nu-1}
			\pas{
				\frac{
					\cos\of{2\pi\frac{k+t}{\nu}}
				}{
					\cos\of{\frac{2\pi\pas{t-1}}{\nu}}
				}
			}^N
	}\d\,t\notag\\
&= -\frac{\log 2}{\nu}-
	\int_{-\frac{1}{2\nu}}^{\frac{1}{2\nu}}
		\log\of{
			\cos\of{2\pi t}
		}\d t.\label{eq:asym2-forrester-even}
\end{align}

For odd $\nu$, the dominating term is
\bit
\item $\cos\of{\frac{2\pi t}{\nu}}$
	for $0\leq t\leq \frac14$,
\item $\cos\of{\frac{2\pi\pas{t+\frac{\nu-1}{2}}}{\nu}}=-\cos\of{\frac{2\pi\pas{t-\frac12}}{\nu}}$
	for $\frac14\leq t\leq \frac34$,
\item $\cos\of{\frac{2\pi\pas{t-1}}{\nu}}$
	for $\frac32\leq t\leq 1$,
\eit
whence by the same simple computation as above, we obtain 
\begin{equation}
\label{eq:asym2-forrester-odd}
F_\nu =-\frac{\log 2}{\nu}-
	2\int_{-\frac{1}{4\nu}}^{\frac{1}{4\nu}}
		\log\of{
			\cos\of{2\pi t}
		}\d t.
\end{equation}

%% file: applics2.tex
\subsection{Gessel and Krattenthaler's formula}
Gessel and Krattentha\-ler \cite{gessel-kratt:cylindric-partitions}
consider nonintersecting paths in the lattice
$\Z\times\Z$, too. However, their lattice paths consist of horizontal and
vertical steps, which essentially is
equivalent to the situation of lattice paths consisting
of diagonal steps in
the even--numbered sublattice
$2\Z\times2\Z$ 
(see Remark~\ref{rem:evenonly}).

\begin{figure}
\caption{Illustration of nonintersecting lattice paths with nonintersecting translate. The shift vector $\mathbf S$ is
indicated by the dotted arrow. }
\label{fig:kratt-paths}
\centerline{\epsfxsize11cm\epsffile{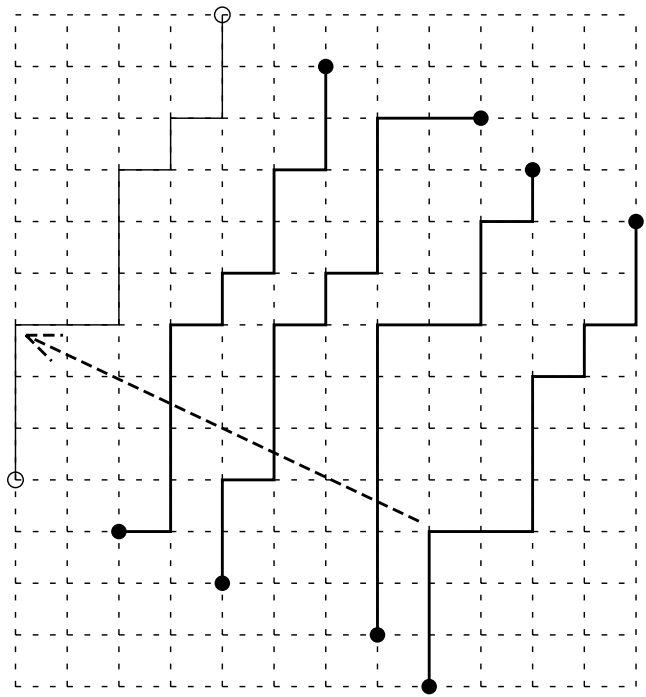}}
\end{figure}

More precisely,
they consider lattice paths which are nonintersecting and, in addition,
are also nonintersecting with respect to ``shifted
copies'' of lattice paths; i.e., copies of the original paths which are
translated by a fixed (non--vertical and non--horizontal) shift vector
$\mathbf S$. See \figref{fig:kratt-paths} for an illustration, where
the translation $\mathbf S$ is indicated by a dotted arrow.

They give a quite general formula 
\cite[Proposition 1, Equation (3.5)]{gessel-kratt:cylindric-partitions}
for the generating function of such nonintersecting families,
in the form of 
a multi--sum of certain determinants.

A special case of this formula 
\bit
\item with shift vector ${\mathbf S}=\pas{-m,m}$,
\item with a certain choice of edge--weights,
\item and with starting points and end points arranged on
	downward--sloping lines
\eit
basically appears as refinement of our formula \eqref{eq:det-gfxy}, in the sense
that now we are only interested in terms with {\em fixed\/} sum ${\sum o = c}$
of offsets of endpoints. So, this amounts to extracting the coefficient
of $y^{c}$ in the expansion of \eqref{eq:det-gfxy}.
The advantage of our formula \eqref{eq:det-gfxy} is that it consists of a
{\em single\/} determinant. Moreover, it does not appear to be
easy to obtain it
by appropriately summing up
Gessel and Krattenthaler's formula.

In any case, the most natural way of understanding \eqref{eq:det-gfxy}
is the {\em direct\/} application of the \lgt\ method.

\subsection*{Acknowledgements}
I thank Christian Krattenthaler for many helpful hints
and discussions.